\documentclass[amstex,12pt, amssymb]{article}

\usepackage{mathtext}
\usepackage[cp1251]{inputenc}
\usepackage[T2A]{fontenc}
\usepackage[dvips]{graphicx}
\usepackage{amsmath}
\usepackage{amssymb}
\usepackage{amsxtra}
\usepackage{latexsym}
\usepackage{ifthen}

\newcounter{lemma}[section]

\newcounter{corollary}[section]

\newcounter{remark}[section]

\newcounter{theorem}[section]

\newcounter{proposition}[section]

\newcounter{example}

\numberwithin{equation}{section}

\begin{document}

\markboth{E.~SEVOST'YANOV, D.~ROMASH, N.~ILKEVYCH}{\centerline{ON
LOWER DISTANCE ESTIMATES ...}}

\def\cc{\setcounter{equation}{0}
\setcounter{figure}{0}\setcounter{table}{0}}

\overfullrule=0pt


\author{EVGENY SEVOST'YANOV, DENYS ROMASH, NATALIYA ILKEVYCH}

\title{
{\bf ON LOWER DISTANCE ESTIMATES OF MAPPINGS IN METRIC SPACES}}

\date{\today}
\maketitle

\begin{abstract}
We consider mappings satisfying an upper bound for the distortion of
families of curves. We establish lower bounds for the distortion of
distances under such mappings. As applications, we obtain theorems
on the discreteness of the limit mapping of a sequence of mappings
converging locally uniformly. We separately consider cases when
mappings are defined in Euclidean $n$-dimensional space and in a
metric space.
\end{abstract}


\section{Introduction}

Lower estimates for distance distortion are an important factor in
providing a description of the local and boundary behavior of many
classes of mappings. In particular, such estimates are known for
quasiconformal mappings, as well as for some more general classes of
mappings, see, e.g., \cite[Theorem~4.7.III]{Ri} and
\cite[Theorem~4.4]{MRSY}. Under studying quasiconformal mappings,
V\"{a}is\"{a}l\"{a} used a lower bound for a certain quantity that
is not directly related to the mapping, but is indirectly related to
the distortion of the distance, see~\cite[Theorem~12.7]{Va}. The
second co-author in a joint work with Ryazanov used a slightly
different lower bound for the capacity, and this gave a result about
the discreteness of the limit map, see the proof of Lemma~4.1
in~\cite{RS}. In particular, the following result holds,
see~\cite[Theorem~4.4]{MRSY}.

\medskip
{\bf Theorem~A.} {\it\, Let $f:{\Bbb B}^n \rightarrow
\overline{{\Bbb R}^n}$ be a $Q$-homeomorphism with $Q\in L^1({\Bbb
B}^n),$ $f(0)=0,$ $h(\overline{{\Bbb R}^n}\setminus f({\Bbb B}^n))
\geqslant \delta
>0$ and $h(f(x_0),f(0)) \geqslant \delta$ for some $x_0 \in {\Bbb B}^n$. Then,
for all $|x|<r=$ $\min(|x_0|/2, 1-|x_0|),$
$$|f(x)|\geqslant \psi(|x|)$$
where $\psi(t)$ is a strictly increasing function with $\psi(0)=0$
which depends only on the $L^1$-norm of $Q$ in ${\Bbb B}^n,$ $n$ and
$\delta$.}

\medskip
Here we use the notion of a $Q$-homeomorphism, which may be found
in~\cite{MRSY}, and, in addition, $h$ is a chordal metric in
$\overline{{\Bbb R}^n}$ (see, e.g., \cite[Section~12]{Va}). A
somewhat more general definition of a ring $Q$-map is given below
for the case of metric spaces, which is also valid for the case of
Euclidean $n$-dimensional space. Our definition of a ring
$Q$-homeomorphism differs slightly from the traditional one, which
is used in~\cite{MRSY}. Note that one of the main goals of this
manuscript is to obtain a more general version of Theorem~$A_1$,
which will be separately formulated and established both in
Euclidean space and in more general metric spaces. Due to this, let
us consider the following definitions.

\medskip
In what follows, $(X, d, \mu)$ and $\left(X^{\,\prime},
d^{\,\prime}, \mu^{\,\prime}\right)$ are metric spaces with metrics
$d$ and $d^{\,\prime}$ and locally finite Borel measures $\mu$ and
$\mu^{\,\prime},$ correspondingly. Let $G$ and $G^{\,\prime}$ be
domains with finite Hausdorff dimensions $\alpha$ and
$\alpha^{\,\prime}\geqslant 2$ in spaces $(X,d,\mu)$ and
$(X^{\,\prime},d^{\,\prime}, \mu^{\,\prime}),$ respectively. For
$x_0\in X$ and $r>0,$ $B(x_0, r)$ and $S(x_0, r)$ denote the ball
$\{x\in X: d(x, x_0)<r\}$ and the sphere $\{x\in X: d(x, x_0)=r\}$,
correspondingly. In what follows,
$${\rm diam}\,E:=\sup\limits_{x, y\in E}d(x, y)\,.$$
\medskip Given $0<r_1<r_2<\infty,$ denote $A=A(x_0, r_1, r_2)=\{x\in
X: r_1<d(x, x_0)<r_2\}.$ Let $p\geqslant$ and $q\geqslant 1,$ and
let $Q:G\rightarrow[0,\infty]$ be a measurable function. Similarly
to \cite[Ch.~7]{MRSY}, a mapping $f:G\rightarrow G^{\,\prime}$ is
called a {\it ring $Q$-mapping at a point $x_0\in \overline{G}$ with
respect to $(p, q)$-moduli}, if for some and for all $0<r_1<
r_2<r_0:={\rm diam\,}G$ the inequality
\begin{equation}\label{eq1C_2}
M_p(f(\Gamma(S(x_0, r_1), S(x_0, r_2), A(x_0, r_1,
r_2))))\leqslant\int\limits_{A(x_0, r_1, r_2)\cap
G}Q(x)\cdot\eta^q(d(x, x_0))\,d\mu(x)
\end{equation}
holds for any measurable function
$\eta:(r_1, r_2)\rightarrow [0, \infty]$ with
\begin{equation}\label{eq8G}
\int\limits_{r_1}^{r_2}\eta(r)\,dr\geqslant 1\ .
\end{equation}
We say that $f$ is $f:G\rightarrow G^{\,\prime}$ is called a {\it
ring $Q$-mapping at a point $x_0\in \overline{G},$ if the latter is
true for $p=\alpha^{\,\prime}$ and $q=\alpha.$}

\medskip
The definition given above, of course, carries over to the case of
the space $X=X^{\,\prime}={\Bbb R}^n,$ $n\geqslant 2.$ In this case,
we set $d(x, y)=d^{\,\prime}(x, y)=|x-y|,$ and, in addition
$\mu=\mu^{\,\prime}=m,$ where $m$ is a Lebesgue measure. We recall
yet another important definition. Due to~\cite{GM}, a domain $D$ in
${\Bbb R}^n$ is called a {\it quasiextremal distance domain} (a
$QED$-{\it domain for short)} if
\begin{equation}\label{eq4***}
M(\Gamma(E, F, {\Bbb R}^n))\leqslant  A\cdot M(\Gamma(E, F, D))
\end{equation}
for some finite number $A\geqslant 1$ and all continua $E$ and $F$
in $D$. In the same way, one can define quasiextremal distance
domains in an arbitrary metric measure space.

\medskip
Given a compact set $K$ in a domain $D$ we set $d(K, \partial
D)=\inf\limits_{x\in K, y\in \partial D}d(x, y).$ If $\partial
D=\varnothing,$ we set $d(K, \partial D)=\infty.$

\medskip
Given a domain $D$ in ${\Bbb R}^n,$ $n\geqslant 2,$ a Lebesgue
measurable function $Q:D\rightarrow [0, \infty],$ a compact set
$K\subset D$ and numbers $A, \delta>0$ denote by $\frak{F}^{A,
\delta}_{K, Q}(D)$ a family of all mappings $f:D\rightarrow {\Bbb
R}^n$ satisfying the relations~(\ref{eq1C_2})--(\ref{eq8G}) at any
point $x_0\in D$ with $d(x, y)=d^{\,\prime}(x, y)=|x-y|$ and
$\mu=\mu^{\,\prime}=m,$ where $m$ is the Lebesgue measure, such that
$D_f:=f(D)$ is a $QED$-domain with $A$ in~(\ref{eq4***}) and, in
addition, $d(f(K), \partial D_f)\geqslant \delta.$ The following
result holds.

\medskip
\begin{theorem}\label{th1}
{\it\,If $Q\in L^1(D),$ then there exist a constant $C>0$ depending
only on $n$ such that
$$|f(x)-f(y)|\geqslant \frac{\delta}{2}\exp\left\{-\frac{\Vert Q\Vert_1A}{C|x-y|^n}\right\}$$
for any $x, y\in K$ and any $f\in \frak{F}^{A, \delta}_{K, Q}(D).$ }
\end{theorem}

\medskip
Theorem~\ref{th1} admits an analogue in metric spaces, which we will
now formulate.

\medskip
Let $X$ be a metric space. We say that the condition of the {\it
complete divergence of paths} is satisfied in $D\subset X$ if for
any different points $y_1$ and $y_2 \in D$ there are some $w_1,$
$w_2\in
\partial D$ and paths $\alpha_2:(-2, -1]\rightarrow
D^{\,\prime},$ $\alpha_1:[1, 2)\rightarrow D$ such that 1)
$\alpha_1$ and $\alpha_2$ are subpaths of some geodesic path
$\alpha: [- 2, 2] \rightarrow X,$ that is, $\alpha_2:=\alpha|_{(-2,
-1]}$ and $\alpha_1:=\alpha|_{[1, 2)};$ 2) 2) the geodesic path
$\alpha$ joins the points $w_2,$ $y_2,$ $y_1$ and $w_1$ such that
$\alpha(-2)=w_2,$ $\alpha(-1)=y_2,$ $\alpha(1)=y_1,$
$\alpha(2)=w_2.$

\medskip
Note that the condition of the complete divergence of the paths is
satisfied for an arbitrary bounded domain $D^{\,\prime}$ of the
Euclidean space ${\Bbb R}^n$, since as paths $\alpha_1$ and
$\alpha_2$ we may take line segments starting at points $y_1$ and
$y_2$ and directed to opposite sides of each other. In this case,
the points $w_1$ and $w_2$ are automatically detected due to the
boundless of $D^{\,\prime}$ (see, e.g., \cite[Proof of
Theorem~1.5]{SevSkv$_1$}.

\medskip
Let $(X, \mu)$ be a metric space with measure $\mu$ and of Hausdorff
dimension $n.$ For each real number $n\geqslant 1,$ we define {\it
the Loewner function} $\phi_n:(0, \infty)\rightarrow [0, \infty)$ on
$X$ as
\begin{equation}\label{eq2H}
\phi_n(t)=\inf\{M_n(\Gamma(E, F, X)): \Delta(E, F)\leqslant t\}\,,
\end{equation}
where the infimum is taken over all disjoint nondegenerate continua
$E$ and $F$ in $X$ and
$$\Delta(E, F):=\frac{{\rm dist}\,(E,
F)}{\min\{{\rm diam\,}E, {\rm diam\,}F\}}\,.$$
A pathwise connected metric measure space $(X, \mu)$ is said to be a
{\it Loewner space} of exponent $n,$ or an $n$-Loewner space, if the
Loewner function $\phi_n(t)$ is positive for all $t> 0$ (see
\cite[Section~2.5]{MRSY} or \cite[Ch.~8]{He}). Observe that ${\Bbb
R}^n$ and ${\Bbb B}^n\subset {\Bbb R}^n$ are Loewner spaces (see
\cite[Theorem~8.2 and Example~8.24(a)]{He}).

\medskip
Let $\Phi_n(t)$ be a Loewner function in~(\ref{eq2H}), corresponding
to the space $U.$ Then, by virtue of~\cite[Theorem~8.23]{He}, there
is $\delta_0> 0$ and some constant $C> 0$ such that
\begin{equation}\label{eq3C}
\Phi_n(t)\geqslant C\log\frac{1}{t}\quad \forall\,\,t>0:
|t|<\delta_0\,.
\end{equation}

\medskip
Given a domain $D$ in $X,$ $n\geqslant 2,$ a measurable function
$Q:D\rightarrow [0, \infty],$ a compact set $K\subset D$ and numbers
$A, \delta>0$ denote by $\frak{F}^{A, \delta}_{K, Q}(D)$ a family of
all mappings $f:D\rightarrow X^{\,\prime}$ satisfying the
relations~(\ref{eq1C_2})--(\ref{eq8G}) at any point $x_0\in D,$ such
that $D_f:=f(D)$ is a compact $QED$-subdomain of $X^{\,\prime}$ with
$A$ in~(\ref{eq4***}) and, in addition, $d^{\,\prime}(f(K),
\partial D_f)\geqslant \delta.$ The following result holds.

\medskip
\begin{theorem}\label{th2}
{\it\, Let $(X, d, \mu)$ and $\left(X^{\,\prime}, d^{\,\prime},
\mu^{\,\prime}\right)$ are metric spaces with metrics $d$ and
$d^{\,\prime}$ and locally finite Borel measures $\mu$ and
$\mu^{\,\prime},$ correspondingly. Let $X^{\,\prime}$ be a Loewner
space. If $Q\in L^1(D),$ then there exist a constant $C>0$ such that
$$d^{\,\prime}(f(x), f(y))\geqslant
\frac{\delta}{2}\exp\left\{-\frac{\Vert Q\Vert_1A}{Cd^n(x,
y)}\right\}$$
for any $x, y\in K$ and any $f\in \frak{F}^{A, \delta}_{K, Q}(D).$ }
\end{theorem}

\medskip
Lower bounds for the distance distortion play an important role in
establishing the discreteness of the limit mapping, as well as in
proving its homeomorphism. In particular, we have the following,
see~\cite[Corollary~21.3]{Va}, \cite[Theorem~4.2]{RS}.

\medskip
{\bf Theorem~B.} {\it If $f_j: D\rightarrow D_j$ is a sequence of
$K$-qc mappings which converge $c$-uniformly to a mapping
$f:D\rightarrow \overline{{\Bbb R}^n},$ then $f$ is either a
homeomorphism onto a domain $D^{\,\prime}$ or a constant.}

\medskip
{\bf Theorem~C.} {\it Let $D$ be a domain in ${\Bbb R}^n,$
$n\geqslant 2,$ and let $Q:D \rightarrow (0,\infty)$ be a measurable
function such that
$$
\int\limits_{0}^{\varepsilon(x_0)}\frac{dr}{rq_{x_0}^{\frac{1}{n-1}}(r)}=\infty\quad
\quad \forall\,x_0\in D$$
for a positive $\varepsilon(x_0)< {\rm dist}\, (x_0,
\partial D)$ where $q_{x_0}(r)$ denotes the average of
$Q(x)$ over the sphere $|x-x_0|=r.$ Suppose that $f_m,$
$m=1,2,\ldots,$ is a sequence of ring $Q$-ho\-me\-o\-mor\-phisms
from $D$ into ${\Bbb R}^n$ converging locally uniformly to a mapping
$f.$ Then the mapping $f$ is either a constant in $\overline{{\Bbb
R}^n}$ or a homeomorphism into ${\Bbb R}^n.$ }

\medskip
With respect to the present paper, the following results are
consequence of Theorems~\ref{th1} and~\ref{th2}.

\medskip
\begin{theorem}\label{th4}
{\it\, Let $D$ be a domain in ${\Bbb R}^n,$ $n\geqslant 2,$ and let
$f_m:D\rightarrow X^{\,\prime},$ $m=1,2,\ldots, $ be a sequence of
homeomorphisms in the class $\frak{F}^{A, \delta}_{K, Q}(D)$
converging to some a mapping $f:D\rightarrow {\Bbb R}^n$ locally
uniformly. If $K=\overline{G}$ and $G$ is a compact subdomain of
$D,$ and $Q\in L^1(D),$ then $f$ is a homeomorphism in $G.$ }
\end{theorem}

\medskip
\begin{theorem}\label{th3}
{\it\, Let $(X, d, \mu)$ and $\left(X^{\,\prime}, d^{\,\prime},
\mu^{\,\prime}\right)$ are metric spaces with metrics $d$ and
$d^{\,\prime}$ and locally finite Borel measures $\mu$ and
$\mu^{\,\prime},$ correspondingly. Let $D$ be a domain in $X$ and
let $f_n:D\rightarrow X^{\,\prime},$ $n=1,2,\ldots, $ be a sequence
of homeomorphisms of the class $\frak{F}^{A, \delta}_{K, Q}(D)$
which is converge locally uniformly to some mapping $f:D\rightarrow
X^{\,\prime}.$ Let $X^{\,\prime}$ be a Loewner space and let
$K=\overline{G}$ and $G$ is a compact subdomain of $D.$ If $Q\in
L^1(D),$ then $f$ is a discrete in $G.$

If in addition, $X^{\,\prime}$ is locally path connected, then $f$
is open. Besides that, if all balls in $X^{\,\prime}$ are connected,
and all closed balls in $X$ are compact, then $f$ is a
homeomorphism.}
\end{theorem}

\medskip
\begin{remark}
It should be borne in mind that an analogue of Theorem~\ref{th4} was
obtained earlier for similar classes of mappings, see Theorem~C.
However, here more stringent conditions on the function $Q$ are
used. Another analogue is Theorem 7.7 in~\cite{MRSY}, where the
condition of local integrability of $Q$ is used. However, in this
case we are dealing with the so-called strong ring
$Q$-homeomorphisms, the definition of which assumes the fulfillment
of a more stringent condition compared
to~(\ref{eq1C_2})--(\ref{eq8G}). Due to the above, the result
obtained in Theorem~\ref{th4} is new even for Euclidean space. As
for more general metric spaces, we do not know of any analogues of
Theorems~\ref{th4} or~\ref{th3}.
\end{remark}

\section{Proof of Theorem~\ref{th1}}

The following statement holds, see, e.g.,
\cite[Theorem~1.I.5.46]{Ku$_2$}).

\medskip
\begin{proposition}\label{pr2}
{\it\, Let $A$ be a set in a topological space $X.$ If the set $C$
is connected and $C\cap A\ne \varnothing\ne C\setminus A,$ then
$C\cap
\partial A\ne\varnothing.$}
\end{proposition}

\begin{figure}[h]
\centerline{\includegraphics[scale=0.4]{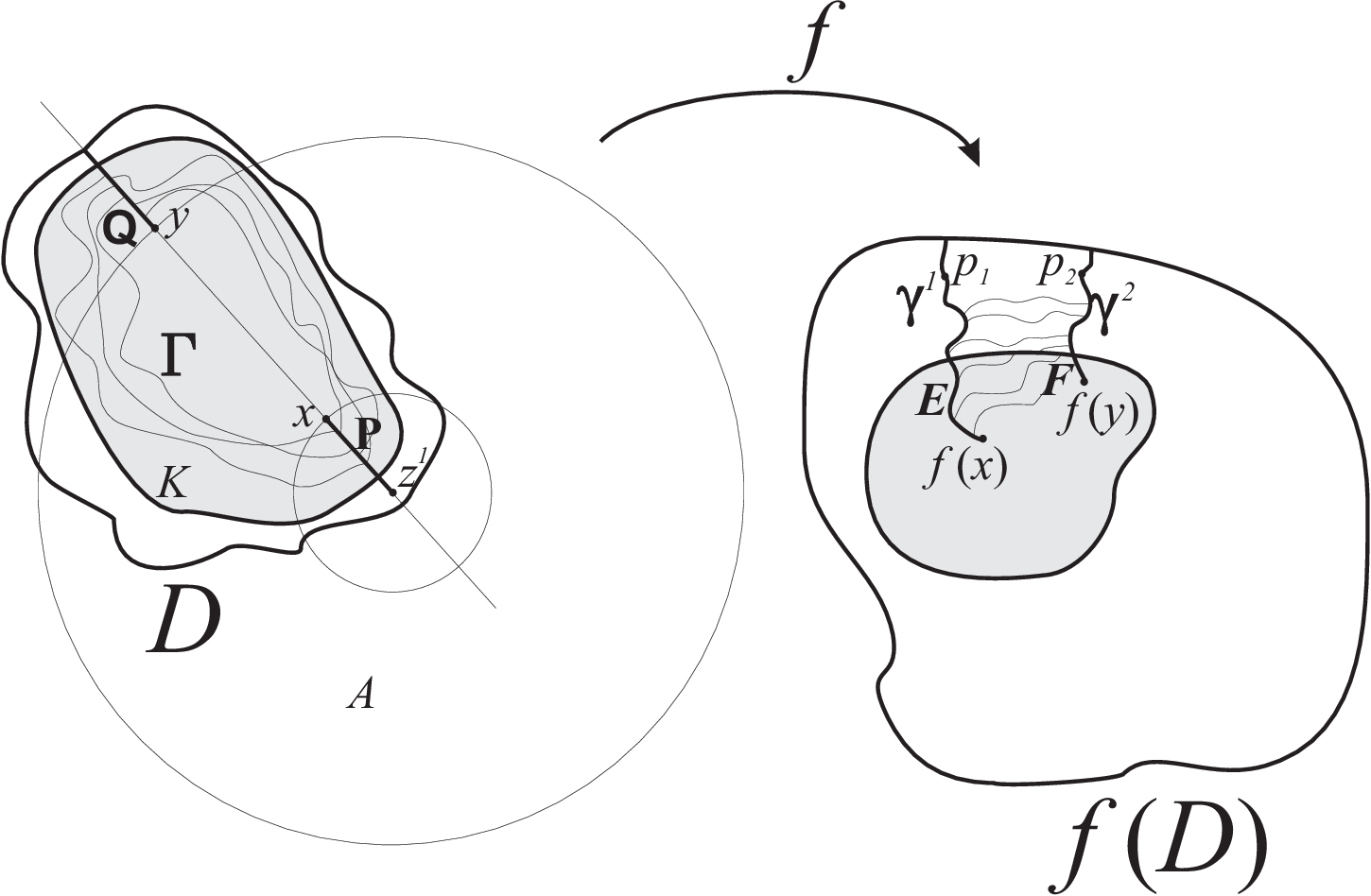}} \caption{To
proof of Theorem~\ref{th2}}\label{fig2}
\end{figure}

\medskip
{\it Proof of Theorem~\ref{th1}.} We make extensive use of our
approaches previously applied to the study of inverse classes of
mappings with respect to $\frak{F}^{A, \delta}_{K, Q}(D),$ see
e.g.~\cite{SevSkv$_1$}--\cite{SevSkv$_2$} and~\cite{SSD}. We fix $x,
y\in K\subset D$ and $f\in \frak{F}^{A, \delta}_{K, Q}(D).$ Draw
through the points $x$ and $y$ a straight line $r=r(t)=x+(x-y)t,$
$-\infty<t<\infty$ (see Figure~\ref{fig2}).
Let $t_1=\sup\limits_{t\geqslant 1: r(t)\in D}t$ and
$t_2:=\inf\limits_{t\leqslant 0: r(t)\in D}t.$ If $r(t)\in D$ for
any $t\geqslant 1,$ we set $t_1=+\infty.$ Similarly, if $r(t)\in D$
for any $t\leqslant 0,$ we set $t_2=-\infty.$ Let $P:=r|_{[1, t_1)}$
and $Q:=r|_{(t_2, 0]}.$ Let $\gamma^1=f(P)$ and $\gamma^2=f(Q).$ Due
to Proposition~13.5 in \cite{MRSY} homeomorphisms preserve the
boundary of a domain. Thus, $\gamma^1(t)\rightarrow \partial f(D)$
as $t\rightarrow t_1-0$ and $\gamma^2(t)\rightarrow \partial f(D)$
as $t\rightarrow t_2+0.$ Let $z^1_m$ be a sequence in $|\gamma^1|,$
$m=1,2,\ldots,$ such that $z^1_m\rightarrow z_0\in \partial f(D).$
If $z_0\ne \infty,$ we may assume that $d(z^1_m, z_0)<1/m.$ Then by
the triangle inequality and due to the definition of the class
$\frak{F}^{A, \delta}_{K, Q}(D)$ we obtain that
$$d(f(x), z^1_m)\geqslant
d(f(x), z_0)-d(z_0, z^1_m)\geqslant \delta-\frac{1}{m}\geqslant
\delta/2\,,\qquad m\geqslant m_0\,.$$
Thus, setting $p_1:=z^1_{m_0}$ we obtain that $d(f(x), p_1)\geqslant
\delta/2.$ If $z_0=\infty,$ the existence of the point $p_1$ on
$|\gamma^1|$ such that $d(f(x), p_1)\geqslant \delta/2$ is obvious.
Similarly, there is a point $p_2$ on $|\gamma^2|$ such that $d(f(x),
p_2)\geqslant \delta/2.$ Let $p_1=\gamma^1(\kappa_1)$ and
$p_1=\gamma^1(\kappa_2)$ for some $\kappa_1\geqslant 1$ and
$\kappa_2\leqslant 0.$ Then $E:=\gamma^1|_{[1, \kappa_1]}$ and
$F:=\gamma^2|_{[\kappa_2, 0]}$ are continua in $f(D)$ such that
$d(E)\geqslant\delta/2$ and $d(F)\geqslant\delta/2.$

\medskip
Next, we should use the fact that the $n$-dimensional Euclidean
space is a Loewner space, see, for example, \cite[Theorem~8.2]{He}.
We will also use the fact that $f(D)$ is a $QED$-domain with the
multiplicative constant $A$ in~(\ref{eq4***}). Observe that,
\begin{equation}\label{eq3D}
\Delta(E, F):=\frac{{\rm dist}\,(E, F)}{\min\{{\rm diam\,}E, {\rm
diam\,}F\}}\leqslant \frac{2|f(x)-f(y)|}{\delta}\,,
\end{equation}
because ${\rm dist}\,(E, F)\leqslant d(f(x), f(y)).$
Then, by the definition of the Loewner function $\Phi_n(t)$
in~(\ref{eq2H}) we obtain that
\begin{equation}\label{eq4C}
\Phi_n\left(\frac{2|f(x)-f(y)|}{\delta}\right)\leqslant M(\Gamma(E,
F, {\Bbb R}^n))\leqslant A\cdot M(\Gamma(E, F, f(D))) \,.
\end{equation}
By~(\ref{eq3C})
\begin{equation}\label{eq5C}
C\cdot\log \frac{\delta}{2|f(x)-f(y)|}\leqslant
\Phi_n\left(\frac{2|f(x)-f(y)|}{\delta}\right)\,.
\end{equation}
Now, by~(\ref{eq4C}) and~(\ref{eq5C}) we obtain that
\begin{equation}\label{eq6C}
C\cdot\log \frac{\delta}{2|f(x)-f(y)|}\leqslant A\cdot M(\Gamma(E,
F, f(D)))\,.
\end{equation}
Now let us prove the upper bound for $M(\Gamma(E, F, f(D))).$
Indeed, by the definition
\begin{equation}\label{eq2}
\Gamma(E, F, f(D))=f(\Gamma(P_1, Q_1, D))\,,
\end{equation}
where $P_1:=f^{\,-1}(E)$ and $Q_1:=f^{\,-1}(F).$ In addition, let us
to show that
\begin{equation}\label{eq1A}
\Gamma(P_1, Q_1, D)>\Gamma(S(z^1, \varepsilon^1), S(z^1,
\varepsilon^2), A(z^1, \varepsilon^1, \varepsilon^2))\,,
\end{equation}
where $z^1:=f^{\,-1}(p_1),$ $\varepsilon^1:=|z^1-x|$ and
$\varepsilon^2:=|z^1-y|.$ Observe that
$$|y-x|+\varepsilon^1=$$
\begin{equation}\label{eq5B}
=|y-x|+|x-z^1|= |z^1-y|=\varepsilon^2\,,
\end{equation}
and, thus, $\varepsilon^1<\varepsilon^2.$ Let $\gamma\in\Gamma(P_1,
Q_1, D).$ Then $\gamma:[0, 1]\rightarrow {\Bbb R}^n,$ $\gamma(0)\in
P_1,$ $\gamma(1)\in Q_1$ and $\gamma(s)\in D$ for $0<s<1.$ Recall
that,
$$z^1=f(y)+(f(x)-f(y))\,\kappa_1$$
for $\kappa_1\geqslant 1.$ Since $\gamma(0)\in P_1,$ there is
$1\leqslant t\leqslant \kappa_1$ such that
$\gamma(0)=f(y)+(f(x)-f(y))t.$ Thus,
$$|\gamma(0)-z^1|=|(x-y)(\kappa_1-t)|$$
\begin{equation}\label{eq2B}
\leqslant |x-y|(\kappa_1-1)=|(x-y)\kappa_1+y-x|\end{equation}
$$=|x-z^1|=\varepsilon^1\,.$$
On the other hand, since $\gamma(1)\in Q_1,$ there is $p\leqslant 0$
such that
$$\gamma(1)=y+(x-y)p\,.$$ In this case, we obtain that
$$
|\gamma(1)-z^1|=|(x-y)(\kappa_1-p)|$$
\begin{equation}\label{eq3A}
\geqslant |(x-y)\kappa_1|=|(x-y)\kappa_1 +y-y|
\end{equation}
$$=|y-z^1|=\varepsilon^2\,.$$
Since $\varepsilon^1<\varepsilon^2,$ due to~(\ref{eq3A}) we obtain
that
\begin{equation}\label{eq4C}
|\gamma(1)-z^1|>\varepsilon^1\,.
\end{equation}
It follows from~(\ref{eq2B}) and (\ref{eq4C}) that $|\gamma|\cap
\overline{B(z^1, \varepsilon^1)}\ne\varnothing\ne (D\setminus
\overline{B(z^1, \varepsilon^1)})\cap|\gamma|.$ In this case, by
Proposition~\ref{pr2} there is $\zeta_1\in (0, 1)$ such
that~$\gamma(\zeta_1)\in S(z^1, \varepsilon^1).$ We may assume that
$\gamma(t)\not\in B(z^1, \varepsilon^1)$ for $t\in (\zeta_1, 1).$
Put $\alpha^1:=\gamma|_{[\zeta_1, 1]}.$

\medskip
On the other hand, since $\varepsilon^1<\varepsilon^2$ and
$\gamma(\zeta_1)\in S(z^1, \varepsilon^1),$ we obtain that
$|\alpha^1|\cap B(z^1, \varepsilon^2)\ne\varnothing.$
By~(\ref{eq3A}) we obtain that $(D\setminus B(z^1,
\varepsilon^2))\cap |\alpha^1|\ne\varnothing.$ Thus, again by
Proposition~\ref{pr2} there is $\zeta_2\in (\zeta_1, 1)$ such that
$\alpha^1(\zeta_2)\in S(z^1, \varepsilon^2).$ We may assume that
$\gamma(t)\in B(z^1, \varepsilon^2)$ for $t\in (\zeta_1, \zeta_2).$
Set $\alpha^2:=\alpha^1|_{[\zeta_1, \zeta_2]}.$ Now,
$\gamma>\alpha^2$ and $\alpha^2\in \Gamma(S(z^1, \varepsilon^1),
S(z^1, \varepsilon^2), A).$ Thus, (\ref{eq1A}) is proved.

\medskip
Now, we set
$$\eta(t)= \left\{
\begin{array}{rr}
\frac{1}{|x-y|}, & t\in [\varepsilon^1, \varepsilon^2],\\
0,  &  t\not\in [\varepsilon^1, \varepsilon^2]\,.
\end{array}
\right. $$
Note that, $\eta$ satisfies the relation~(\ref{eq8G}) for
$r_1=\varepsilon^1$ and $r_2=\varepsilon^2.$ Indeed, it follows
from~(\ref{eq5B}) that
$$r_1-r_2=\varepsilon^2-\varepsilon^1=|y-z^1|-|x-
z^1|=$$$$=|x-y|\,.$$
Then
$\int\limits_{\varepsilon^1}^{\varepsilon^2}\eta(t)\,dt=(1/|x-y|)\cdot
(\varepsilon^2-\varepsilon^1)= 1.$ By the inequalities~(\ref{eq2})
and (\ref{eq1A}), and by the definition of the class $\frak{F}^{A,
\delta}_{K, Q}(D),$ we obtain that
$$M(\Gamma(E, F, f(D)))\leqslant M(f(\Gamma(S(z^1, \varepsilon^1), S(z^1,
\varepsilon^2), A(z^1, \varepsilon^1, \varepsilon^2))))$$
\begin{equation}\label{eq14***}
\leqslant \frac{1}{|x-y|^n}\int\limits_{D} Q(x)\,dm(x)=\frac{\Vert
Q\Vert_1}{{|x-y|}^n}\,.
\end{equation}
By~(\ref{eq6C}) and (\ref{eq14***}), we obtain that
$$\frac{C}{A}\cdot\log \frac{\delta}{2|f(x)-f(y)|}\leqslant
\frac{\Vert Q\Vert_1}{{|x-y|}^n}\,.$$
From the latter ratio, the desired inequality directly
follows.~$\Box$

\section{Proof of Theorem~\ref{th2}}

The proof of Theorem~\ref{th2} differs insignificantly from the
proof of Theorem~\ref{th1}, so we will omit unnecessary details to
the extent that they are necessary and will limit ourselves in this
case to elements of a schematic proof.

\medskip
We fix $x, y\in K\subset D$ and $f\in \frak{F}^{A, \delta}_{K,
Q}(D).$ By hypothesis, for the points $x$ and $y\in D$ there are
$w_1,$ $w_2\in
\partial D$ and paths $\alpha_2:(-2, -1]\rightarrow D,$ $\alpha_1:[1,
2)\rightarrow D$ such that

\medskip
1) $\alpha_1$ and $\alpha_2$ are subpaths of some geodesic path
$\alpha:[-2, 2]\rightarrow X,$ that is, $\alpha_2:=\alpha|_{(-2,
-1]}$ and $\alpha_1:=\alpha|_{[1, 2)};$

2) the geodesic $\alpha$ sequentially joins the points  $w_2,$ $x,$
$y$ and $w_1,$ namely, $\alpha(-2)=w_2,$ $\alpha(-1)=x,$
$\alpha(1)=y,$ $\alpha(2)=w_2.$ Note that, for each $i=1,2$ the sets
$|\alpha_2|:=\{x\in X: \exists\,t\in (-2,-1]: \alpha_2(t)=x\}$ and
$|\alpha_1|:=\{x\in X: \exists\,t\in [1, 2): \alpha_1(t)=x\}$ are
connected as continuous images of the corresponding intervals $(-2,
-1]$ and $[1, 2)$ under the mapping $f.$

\medskip
Let $P:=\alpha|_{[1, 2)}$ and $Q:=\alpha|_{(-2, 1]}.$ Let
$\gamma^1=f(P)$ and $\gamma^2=f(Q).$ Since $D^{\,\prime}_f=f(D)$ has
a compact closure in $X^{\,\prime}$ by the assumption, the cluster
set $C(w_i, f)\ne \varnothing,$ $i=1,2,$ in addition, $C(w_i,
f)\subset \partial f(D)$ (see \cite[Proposition~13.5]{MRSY}). Let
$z^1_m$ be a sequence in $|\gamma^1|,$ $m=1,2,\ldots,$ such that
$z^1_m\rightarrow z_0\in \partial f(D).$ We may assume that
$d^{\,\prime}(z^1_m, z_0)<1/m.$ Then by the triangle inequality and
due to the definition of the class $\frak{F}^{A, \delta}_{K, Q}(D)$
we obtain that
$$d^{\,\prime}(f(x), z^1_m)\geqslant
d^{\,\prime}(f(x), z_0)-d^{\,\prime}(z_0, z^1_m)\geqslant
\delta-\frac{1}{m}\geqslant \delta/2\,,\qquad m\geqslant m_0\,.$$
Thus, setting $p_1:=z^1_{m_0}$ we obtain that $d^{\,\prime}(f(x),
p_1)\geqslant \delta/2.$ Similarly, there is a point $p_2$ on
$|\gamma^2|$ such that $d^{\,\prime}(f(x), p_2)\geqslant \delta/2.$
Let $p_1=\gamma^1(\kappa_1)$ and $p_1=\gamma^1(\kappa_2)$ for some
$\kappa_1\geqslant 1$ and $\kappa_2\leqslant 0.$ Then
$E:=\gamma^1|_{[1, \kappa_1]}$ and $F:=\gamma^2|_{[\kappa_2, 0]}$
are continua in $f(D)$ such that $d^{\,\prime}(E)\geqslant\delta/2$
and $d(F)\geqslant\delta/2.$

\medskip
Recall that, Since $X^{\,\prime}$ is a Loewner space and, besides
that, $f(D)$ is a $QED$-domain with the multiplicative constant $A$
in~(\ref{eq4***}). Observe that,
\begin{equation}\label{eq3E}
\Delta(E, F):=\frac{{\rm dist}\,(E, F)}{\min\{{\rm diam\,}E, {\rm
diam\,}F\}}\leqslant \frac{2d^{\,\prime}(f(x),f(y))}{\delta}\,,
\end{equation}
because ${\rm dist}\,(E, F)\leqslant d^{\,\prime}(f(x), f(y)).$
Then, by the definition of the Loewner function
$\Phi_{\alpha^{\,\prime}}(t)$ in~(\ref{eq2H}) we obtain that
\begin{equation}\label{eq4D}
\Phi_{\alpha^{\,\prime}}\left(\frac{2d^{\,\prime}(f(x),f(y))}{\delta}\right)\leqslant
M(\Gamma(E, F, X^{\,\prime}))\leqslant A\cdot M(\Gamma(E, F, f(D)))
\,.
\end{equation}
By~(\ref{eq3E})
\begin{equation}\label{eq5D}
C\cdot\log \frac{\delta}{2d^{\,\prime}(f(x),f(y))}\leqslant
\Phi_{\alpha^{\,\prime}}\left(\frac{2d^{\,\prime}(f(x),f(y))}{\delta}\right)\,.
\end{equation}
Now, by~(\ref{eq4D}) and~(\ref{eq5D}) we obtain that
\begin{equation}\label{eq6D}
C\cdot\log \frac{\delta}{2d^{\,\prime}(f(x),f(y))}\leqslant A\cdot
M(\Gamma(E, F, f(D)))\,.
\end{equation}
Now let us prove the upper bound for $M(\Gamma(E, F, f(D))).$
Indeed, by the definition
\begin{equation}\label{eq2A}
\Gamma(E, F, f(D))=f(\Gamma(P_1, Q_1, D))\,,
\end{equation}
where $P_1:=f^{\,-1}(E)$ and $Q_1:=f^{\,-1}(F).$ Using the similar
approach which was used under the proof of the
relation~(\ref{eq1A}), we may show that
\begin{equation}\label{eq1B}
\Gamma(P_1, Q_1, D)>\Gamma(S(z^1, \varepsilon^1), S(z^1,
\varepsilon^2), A(z^1, \varepsilon^1, \varepsilon^2))\,,
\end{equation}
where $z^1:=f^{\,-1}(p_1),$ $\varepsilon^1:=d(z^1, x)$ and
$\varepsilon^2:=d(z^1, y).$

\medskip
Now, we set
$$\eta(t)= \left\{
\begin{array}{rr}
\frac{1}{d(x, y)}, & t\in [\varepsilon^1, \varepsilon^2],\\
0,  &  t\not\in [\varepsilon^1, \varepsilon^2]\,.
\end{array}
\right. $$
Note that, $\eta$ satisfies the relation~(\ref{eq8G}) for
$r_1=\varepsilon^1$ and $r_2=\varepsilon^2.$ Indeed,
$$\int\limits_{r_1}^{r_2}\eta(t)\,dt=\frac{r_2-r_1}{d(x, y)}=
\frac{d(z_1, x)-d(z_1, y)}{d(x, y)}=1\,,$$ since all three points
$x,$ $y$ and $z_1$ are sequentially located on one geodesic, which
means that
$$r_2=d(z_1,
x)=d(z_1, x)+d(x, y)=r_1+d(x, y)\,.$$
By the inequalities~(\ref{eq2}) and (\ref{eq1A}), and by the
definition of the class $\frak{F}^{A, \delta}_{K, Q}(D),$ we obtain
that
$$M(\Gamma(E, F, f(D)))\leqslant M(f(\Gamma(S(z^1, \varepsilon^1), S(z^1,
\varepsilon^2), A(z^1, \varepsilon^1, \varepsilon^2))))$$
\begin{equation}\label{eq14****}
\leqslant \frac{1}{d^n(x, y)}\int\limits_{D}
Q(x)\,d\mu(x)=\frac{\Vert Q\Vert_1}{{d}^n(x, y)}\,.
\end{equation}
By~(\ref{eq6D}) and (\ref{eq14***}), we obtain that
$$\frac{C}{A}\cdot\log \frac{\delta}{2d^{\,\prime}(f(x),f(y))}\leqslant
\frac{\Vert Q\Vert_1}{{d}^n(x, y)}\,.$$
From the latter ratio, the desired inequality directly
follows.~$\Box$

\section{Convergence theorems}

{\it Proof of Theorem~\ref{th4}.} Let us firstly prove that $f$ is
discrete in $G.$ Indeed, if $f$ is not discrete, then there is a
point $x_0\in G$ and a sequence $x_k\in G,$ $x_k \ne x_0,$
$k=1,2\ldots,$ such that $x_k \rightarrow x_0$ as $k\rightarrow
\infty$ with $f(x_k)=f(x_0).$ By Theorem~\ref{th1}
$$|f_m(x)-f_m(x_0)|\geqslant \frac{\delta}{2}\exp\left\{-\frac{\Vert Q\Vert_1A}{C|x-x_0|^n}\right\}\,,
\quad m=1,2,\ldots$$
in some neighborhood of $x_0.$ Passing to the limit as
$m\rightarrow\infty$ in the last inequality, we obtain that
$$|f(x)-f(x_0)|\geqslant \frac{\delta}{2}\exp\left\{-\frac{\Vert Q\Vert_1A}{C|x-x_0|^n}\right\}\,.$$
In particular,
$$0=|f(x_k)-f(x_0)|\geqslant
\frac{\delta}{2}\exp\left\{-\frac{\Vert
Q\Vert_1A}{C|x_k-x_0|^n}\right\}>0\,.$$
The above contradiction proves the discreteness of $f.$ Finally, $f$
is a homeomorphism due to~\cite[Theorem~3.1]{RS}.~$\Box$

\medskip
{\it Proof of Theorem~\ref{th4}.} Let us use the methodology used in
proving Theorem~7.7 in~\cite{MRSY}. Just as above, it may be shown
that the limit mapping $f$ is discrete. Indeed, if $f$ is not
discrete, then there is a point $x_0\in G$ and a sequence $x_k\in
G,$ $x_k \ne x_0,$ $k=1,2\ldots,$ such that $x_k \rightarrow x_0$ as
$k\rightarrow \infty$ with $f(x_k)=f(x_0).$ By Theorem~\ref{th2}
$$d^{\,\prime}(f_n(x),f_n(x_0))
\geqslant \frac{\delta}{2}\exp\left\{-\frac{\Vert
Q\Vert_1A}{Cd^n(x,x_0)}\right\} \quad n=1,2,\ldots$$
in some neighborhood of $x_0.$ Passing to the limit as
$n\rightarrow\infty$ in the last inequality, we obtain that
\begin{equation}\label{eq3}
d^{\,\prime}(f(x),f(x_0))\geqslant
\frac{\delta}{2}\exp\left\{-\frac{\Vert
Q\Vert_1A}{Cd^n(x,x_0)}\right\}\,.
\end{equation}
In particular,
$$0=d^{\,\prime}(f(x_k),f(x_0))\geqslant
\frac{\delta}{2}\exp\left\{-\frac{\Vert Q\Vert_1A}{Cd^n(x_k,
x_0)}\right\}>0\,.$$
The above contradiction proves the discreteness of $f.$

\medskip
Let us also prove that $f$ is open with respect to the metric
$d^{\,\prime},$ i.e., given an open set $A\subset G,$ for any
$y_0\in f(A)$ there is $\varepsilon_1>0$ such that $B(y_0,
\varepsilon_1)\subset f(A).$ Let $A$ and $y_0$ be as stated above.
Let $x_0\in A$ and let $\varepsilon_0>0$ be such that $f(x_0)=y_0$
and $B(x_0, \varepsilon_0)\subset A.$ It is sufficient to prove
that, there is $\varepsilon_1>0$ such that $B(y_0,
\varepsilon_1)\subset f(B(x_0, \varepsilon_0)).$ Assume the
contrary, namely, given $k\in {\Bbb N}$ there is $y_k\in B(y_0,
1/k)$ such that $y_k\not\in f(B(x_0, \varepsilon_0)).$

\medskip
Since $B(y_0, 1/k)$ is locally path connected by the assumption,
there are a sequence $V_k$ of path connected neighborhoods of $y_0$
such that $V_k\subset B(y_0, 1/k)$ for any $k\in {\Bbb N}.$ Without
loss of generality we may assume that $y_k\in V_k,$ $k=1,2,\ldots .$
Let $\gamma_k$ be a path joining $y_k$ and $y_0$ in $X^{\,\prime},$
$\gamma_k:[0, 1]\rightarrow V_k,$ $\gamma_k(0)=y_0$ and
$\gamma_k(1)=y_k$ for any $k\in {\Bbb N}.$ Let $t_k=\sup_{t\in [0,
1]:\gamma_k(t)\in f(B(x_0, \varepsilon))}t$ and let
$\gamma^{\,\prime}_k:=\gamma_k|_{[0, t_k)}.$ Set
$\gamma^{*}_k:=f^{\,-1}(\gamma^{\,\prime}_k).$ Observe that,
$\gamma^{*}_k(t)\rightarrow S(x_0, \varepsilon_0)$ (see e.g.
Proposition~13.5 in \cite{MRSY}). Thus, there is a sequence $x_k$ in
$B(x_0, \varepsilon_0)$ such that $x_k\in |\gamma_k|,$ $d(x_k,
x_0)\geqslant \varepsilon_0/2.$ Since $f(x_k)\in V_k\subset B(y_0,
1/k),$ we have that $d^{\,\prime}(f(x_k), y_0)\rightarrow 0$ as
$k\rightarrow\infty.$ On the other hand, by~(\ref{eq3})
$$d^{\,\prime}(f(x_k), y_0)\geqslant
\frac{\delta}{2}\exp\left\{-\frac{\Vert Q\Vert_1A}{Cd^n(x_k,
x_0)}\right\}\geqslant \frac{\delta}{2}\exp\left\{-\frac{\Vert
Q\Vert_1A2^n}{C\varepsilon^n_0}\right\}\,.$$
The obtained contradiction proves the assumption made above. Thus,
$f$ is open, as required.

\medskip
Now we will show that the mapping $f$ is injective under the
additional assumption that any ball in $X^{\,\prime}$ is connected.
Let us reason by contradiction. Let us assume that there are $x_1,
x_2\in G$ such that $f(x_1)\ne f(x_2).$ Set $B_t:=B(x_1, t).$ Due to
the discreteness of $f,$ there is $t_0>0$ such that $x_2\not\in
\overline{B(x_1, t)}$ for any $t\in [0, t_0].$

\medskip
Let $y_m:=f_m(x_1),$ $z_m=f_m(x_2)$ and $C_m:=f_m(B_t).$ Let us to
show that,
\begin{equation}\label{eq1}
B(y_m, d^{\,\prime}(y_m, \partial C_m))\subset C_m\,,\qquad
m=1,2,\ldots\,.
\end{equation}
Assume the contrary. Then is $y\in B(y_m, d^{\,\prime}(y_m, \partial
C_m))\setminus C_m.$ Now, $$B(y_m, d^{\,\prime}(y_m, \partial
C_m))\setminus C_m\ne\varnothing\ne C_m\cap B(y_m, d^{\,\prime}(y_m,
\partial C_m))$$ because $y_m\in C_m,$ as well.
By Proposition~\ref{pr2}, since by the assumption the ball $B(y_m,
d^{\,\prime}(y_m,
\partial C_m))$ is connected we obtain that $B(y_m,
d^{\,\prime}(y_m,
\partial C_m))\cap \partial C_m\ne\varnothing.$ Now, there is $z\in B(y_m,
d^{\,\prime}(y_m,
\partial C_m))\cap \partial C_m.$ Thus,
$$d^{\,\prime}(y_m,
\partial C_m)\leqslant d^{\,\prime}(z, y_m)<d^{\,\prime}(y_m,
\partial C_m)\,.$$
The contradiction obtained above proves~(\ref{eq1}). Now,
by~(\ref{eq1})
\begin{equation}\label{eq2C}
\overline{B(y_m, d^{\,\prime}(y_m, \partial C_m))}\subset
\overline{C_m}\,,\qquad m=1,2,\ldots\,.
\end{equation}

\medskip
Due to~(\ref{eq2C}), $d^{\,\prime}(y_m, \partial
C_m)<d^{\,\prime}(y_m, z_m),$ $m=1,2,\ldots. $ Indeed, in the
contrary case $d^{\,\prime}(y_m, z_m)\leqslant d^{\,\prime}(y_m,
\partial C_m)$ for some $m\in{\Bbb N},$ so that
$z_m\in \overline{B(y_m, d^{\,\prime}(y_m, \partial C_m))}.$ But
now, by~(\ref{eq2C}) $z_m\in \overline{C_m}$ that contradicts the
choice of $z_m=f_m(x_2)$ and that $f_m$ is a homeomorphism.

\medskip
Since $f_m$ is homeomorphic and open in the topology of
$X^{\,\prime},$ $\partial C_m=f_m(\partial B_t).$ In addition, by
the assumption, the closed balls are compacta in $X,$ so that
$\partial C_m=f_m(\partial B_t)$ is compact set. Now, there is
$x_{m,t} \in \partial B_t$ such that
\begin{equation}\label{eq3A}
d^{\,\prime}(y_m, \partial C_m) = d^{\,\prime}(y_m,
f_m(x_{m,t}))\,,\qquad m=1,2,\ldots\,.
\end{equation}
By compactness of $\partial B_t$, for every $t\in (0, t_0],$ there
is $x_t\in \partial B_t$ such that $d^{\,\prime}(x_{m_k,t},
x_t)\rightarrow 0$ as $k\rightarrow \infty$ for some subsequence
$m_k,$ $k=1,2,\ldots .$ Since the locally uniform convergence of
continuous functions in a metric space implies the continuous
convergence (see Theorem 21.X.3 in \cite{Ku$_1$}), we have that
$d^{\,\prime}(f_{m_k}(x_{m_k,t}),f(x_{t})) \rightarrow 0$ as
$k\rightarrow \infty.$ Consequently, by~(\ref{eq3}) and (\ref{eq3A})
we obtain that
$$d^{\,\prime}(f(x_1),f(x_t))\leqslant d^{\,\prime}(f(x_1),f(x_2))$$
for any $t\in(0, t_0].$
However, since $f(x_1)=f(x_2)$ we obtain that $f(x_t)=f(x_1)$ for
every $t\in (0,t_0].$ The latter contradicts to the discreteness of
$f.$ Thus, $f$ is injective. In addition, $f^{\,-1}$ is continuous
in $f(G)$ because $f$ is open by the proving above.~$\Box$

\medskip
{\bf \noindent Evgeny Sevost'yanov} \\
{\bf 1.} Zhytomyr Ivan Franko State University,  \\
40 Velyka Berdychivs'ka Str., 10 008  Zhytomyr, UKRAINE \\
{\bf 2.} Institute of Applied Mathematics and Mechanics\\
of NAS of Ukraine, \\
19 Henerala Batyuka Str., 84 116 Slov'yansk,  UKRAINE\\
esevostyanov2009@gmail.com

\medskip
{\bf \noindent Denys Romash} \\
Zhytomyr Ivan Franko State University,  \\
40 Velyka Berdychivs'ka Str., 10 008  Zhytomyr, UKRAINE \\
dromash@num8erz.eu

\medskip
{\bf \noindent Nataliya Ilkevych} \\
Zhytomyr Ivan Franko State University,  \\
40 Velyka Berdychivska Str., 10 008  Zhytomyr, UKRAINE \\
ilkevych1980@gmail.com

\end{document}